\documentclass[runningheads]{llncs}
\usepackage[T1]{fontenc}
\usepackage{graphicx}
\usepackage{hyperref}
\usepackage{color}

\usepackage{amsmath}
\usepackage{amssymb}
\usepackage{booktabs}
\usepackage[outdir=./]{epstopdf}
\usepackage{caption}
\usepackage{subcaption}

\newcommand{\setR}{\mathbb{R}}

\newcommand{\ignore}[1]{}

\begin{document}
\title{Efficient Line Search Method Based on Regression and Uncertainty Quantification}
\titlerunning{Efficient Line Search Method}
%
\author{Tomislav Prusina \and
Sören Laue}
%
%
\institute{Universität Hamburg, Germany\\
\email{\{tomislav.prusina, soeren.laue\}@uni-hamburg.de}}
\maketitle              
\begin{abstract}
Unconstrained optimization problems are typically solved using iterative methods, which often depend on line search techniques to determine optimal step lengths in each iteration. This paper introduces a novel line search approach. Traditional line search methods, aimed at determining optimal step lengths, often discard valuable data from the search process and focus on refining step length intervals. This paper proposes a more efficient method using Bayesian optimization, which utilizes all available data points, i.e., function values and gradients, to guide the search towards a potential global minimum. This new approach more effectively explores the search space, leading to better solution quality. It is also easy to implement and integrate into existing frameworks. Tested on the challenging CUTEst test set, it demonstrates superior performance compared to existing state-of-the-art methods, solving more problems to optimality with equivalent resource usage.

\keywords{Nonlinear Optimization \and Line Search \and Regression \and Uncertainty Quantification \and Bayesian Optimization.}
\end{abstract}

\section{Introduction}

Optimization plays a central role in solving many real-world problems, from engineering and operations research to training machine learning models~\cite{GiesenJL12,GiesenL19,sra2012}. The selection of parameters and strategies significantly influences the effectiveness of algorithms designed to tackle such optimization problems. A crucial aspect of these algorithms, which are often iterative methods, is the identification of appropriate step lengths, a process known as line search, which is fundamental to the efficiency of the optimization algorithm.

Over the years, many line search methods have been proposed, each to ensure a step length that facilitates considerable progress in each iteration. Often, these methods try to satisfy the strong Wolfe conditions, ensuring a sufficient decrease in the function value and a reduction in the directional derivative. This ensures that the algorithm converges to a point close to a local minimum of the function along the search ray.

Traditionally, line search methods achieve this by maintaining and iteratively refining an interval with upper and lower bounds for the step length. This 'zooming-in' process is designed to pinpoint a suitable step length that guarantees sufficient progress. However, in this process, only the upper and lower bounds are kept and updated, and the information, i.e., function value and gradient information from all other query points visited during this line search, is discarded. Also, other regions of the line search interval might become more interesting, but they are never considered.

This paper addresses this shortfall. We propose a novel approach based on Bayesian optimization that harnesses all available information – the function value and gradient at all points queried so far – to identify a potential global minimum of the function along the line search ray. 
Bayesian optimization provides a principled framework for efficiently exploring the search space by constructing a probabilistic model of the objective function. By iteratively refining this model based on observed function values, Bayesian optimization guides the search toward promising regions, thus facilitating faster convergence and improved solution quality. 

Our approach is simple, lightweight, and can be seamlessly implemented with a few lines of Python code, making it accessible and easily integrated into existing frameworks. We test our approach on the CUTEst test set~\cite{Gould2015} that contains a number of challenging unconstrained optimization problems. On these problems, our method outperforms other state-of-the-art approaches. Specifically, it enables a solver to solve more problems to optimality, while using an equivalent number of function evaluations and iterations.

\section{Related Work}
Many line search methods have been proposed. The seminal work by Armijo~\cite{armijo66} introduced the Armijo rule, one of the first inexact line search methods. Building on this, Wolfe~\cite{wolfe69} formulated the weak and strong Wolfe conditions, crucial for balancing computational efficiency with sufficient descent and curvature criteria. Barzilai and Borwein~\cite{barzilai88} proposed a line search approach using gradient information, which proved effective in large-scale optimization. The work of Moré and Thuente~\cite{MoreT94}presented a line search method based on cubic interpolation, improving accuracy and efficiency in various applications. It is used in the quasi-Newton solver L-BFGS-B~\cite{ZhuBLN97} which is also part of the Scipy library~\cite{2020SciPy-NMeth}. Additionally, Hager and Zhang~\cite{Hager2005} introduced an approach to refine step size estimation, enhancing performance in nonlinear optimization problems.
Wächter and Biegler~\cite{Waechter2005} introduced a filter line search strategy for primal-dual interior-point methods for large-scale nonlinear problems where traditional line search methods may struggle. It is part of the Ipopt solver~\cite{Waechter2006}, which is especially efficient in handling complex and large-scale optimization problems. 
All these line search methods have in common that they narrow down an interval containing a good step length that ensures adequate progress. Only the upper and lower bounds of this interval are maintained and updated, while the data, specifically function value and gradient information from other points assessed in the line search, is disregarded. Additionally, other potentially relevant areas within the line search interval are not explored.

On the other hand, Bayesian optimization has emerged as a prominent method in global optimization, particularly for optimizing expensive-to-evaluate functions. The pioneering work by Mo\v{c}kus~\cite{mockus1972bayes} laid the foundation for this technique in the one-dimensional case, focusing on efficiency in the face of limited evaluation opportunities. The integration of Gaussian processes in Bayesian optimization, as explored by Rasmussen and Williams~\cite{RasmussenW06}, provided a robust framework for modeling the objective function and quantifying uncertainty. Usually, Bayesian optimization models use only function values and not gradient information. An approach that also integrated gradient information into Gaussian processes in the context of Bayesian optimization is described in~\cite{NIPS2017}. However, their algorithm is rather involved, and no code is provided. All these works did not consider line search as an application for Bayesian optimization. Here, we present a line search approach that uses Bayesian optimization and that is easy to implement.

\section{Line Search Method}
Let $f\colon\setR^n\to\setR$ be a continuously differentiable function and  bounded from below. Suppose we want to solve the unconstrained optimization problem
\begin{equation}
\min f(x)
\end{equation}

\begin{figure}
     \centering
     \begin{subfigure}[b]{0.328\textwidth}
         \centering
         \includegraphics[width=\linewidth]{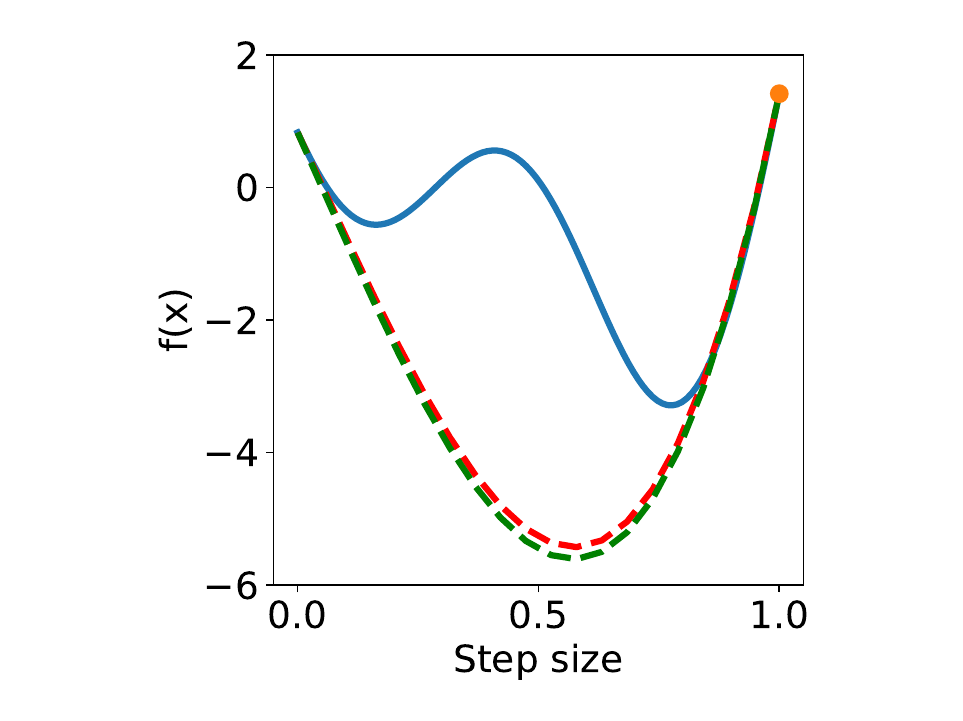}
         \caption{First iteration.}
         \label{fig:plot0}
     \end{subfigure}
\hfill
     \begin{subfigure}[b]{0.328\textwidth}
         \centering
         \includegraphics[width=\linewidth]{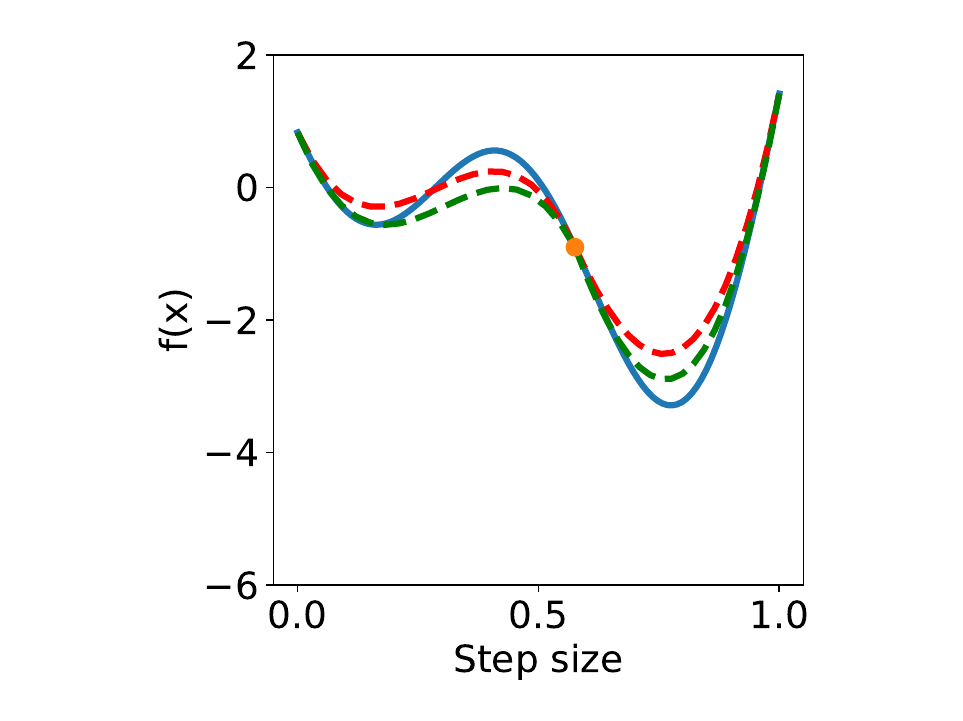}
         \caption{Second iteration.}
         \label{fig:plot1}
     \end{subfigure}
\hfill
     \begin{subfigure}[b]{0.328\textwidth}
         \centering
         \includegraphics[width=\linewidth]{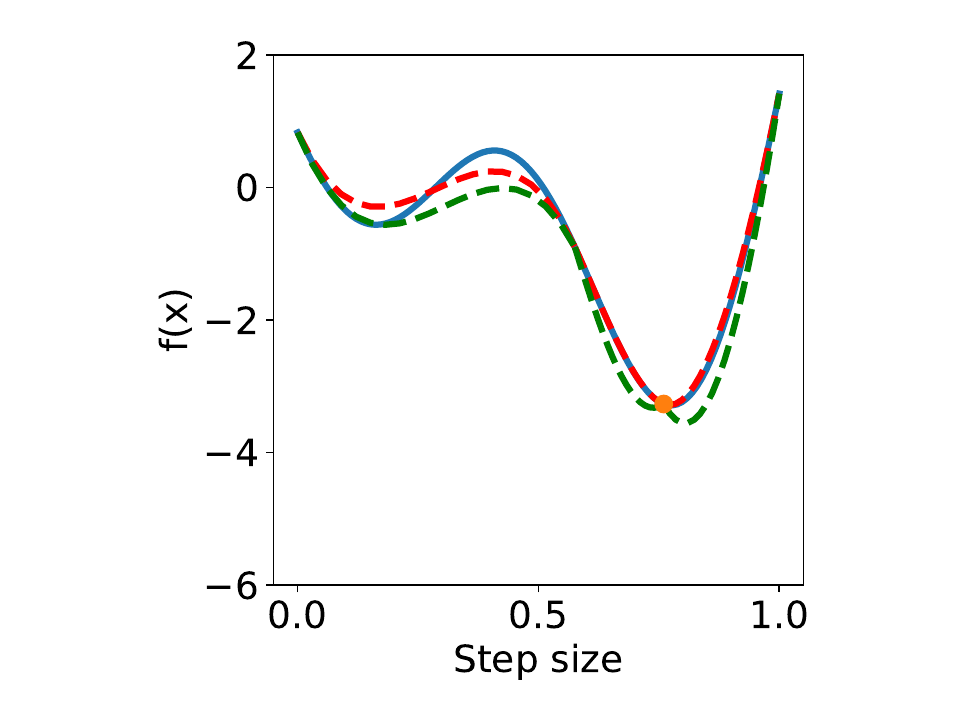}
         \caption{Third iteration.}
         \label{fig:plot2}
     \end{subfigure}
        \caption{The figures show three iterations of our proposed line search method. The blue line is the original function, the dashed red line is the cubic interpolation, and the dashed green line is the uncertainty. The x-axis shows the step length, while the y-axis shows the function value.}
        \label{fig:1}
\end{figure}

Iterative methods usually proceed as follows. Given a starting point $x_0\in\setR^n$, in each iteration the current iterate $x_k$ is updated by
\[
x_{k+1} = x_k +\alpha_k p_k,
\]
where $\alpha_k$ is the step length and $p_k$ is a descent direction, i.e., a direction in which the function $f$ decreases. Mathematically, it can expressed as $p_k^\top \nabla f(x_k) < 0$. An (approximate) line search method tries to find a good estimate of the step length parameter $\alpha_k$ such that sufficient progress is made in the current iteration. In our approach, we use the strong Wolfe conditions to ensure sufficient progress, i.e., the final $\alpha_k$ needs to satisfy
\[
\begin{aligned}
f(x_k+\alpha_k p_k) & \leq f(x_k) + c_1\alpha_k \nabla f(x_k)^\top p_k \\
|\nabla f(x_k+\alpha_k p_k)^\top p_k| & \leq c_2 | \nabla f(x_k)^\top p_k|,
\end{aligned}
\]
where $0 < c_1 < c_2 < 1$ are two positive constants.
The first inequality ensures a sufficient decrease in function value and the second ensures a point that is close to a stationary point along the search ray.  Such a step length can be found by the following optimization problem
\begin{equation}
\min_{\alpha_k\geq 0} f(x_k + \alpha_k p_k),
\label{prob:1}
\end{equation}
i.e., the global minimizer along the search direction satisfies the strong Wolfe conditions. Usually, we do not know the structure of function $f$ but only have access to an oracle that computes function value and gradients for a given query point. Hence, finding an approximate solution to  Problem~\eqref{prob:1} can be done using Bayesian optimization. For this, we create a model of the one-dimensional function along the search ray as follows: We approximate the true function by a cubic spline interpolation and additionally add an uncertainty estimate that depends on the distance to the closest query point and follows a Gaussian distribution. See Fig.~\ref{fig:1} for an illustration. 

Suppose that in the current iteration we have already queried the function for a number of query points. There always exists exactly one solution for a cubic polynomial interpolating between two consecutive query points such that the polynomial matches function values and directional gradients in the given points. The uncertainty estimate for a given point follows a Gaussian distribution that is scaled with the function value of the corresponding cubic polynomial. In each iteration, the minimum of this Bayesian model, i.e., the spline interpolation minus the uncertainty estimate is computed and the result is the next query point. This process is repeated until a step length is found that satisfies the strong Wolfe conditions. In the practical implementation, an upper bound on the number of query points is fixed to account for numerical instabilities and round-off errors. This is common practice in line search methods.

\begin{figure}
     \centering
     \begin{subfigure}[b]{.8\textwidth}
         \centering
         \includegraphics[width=\linewidth]{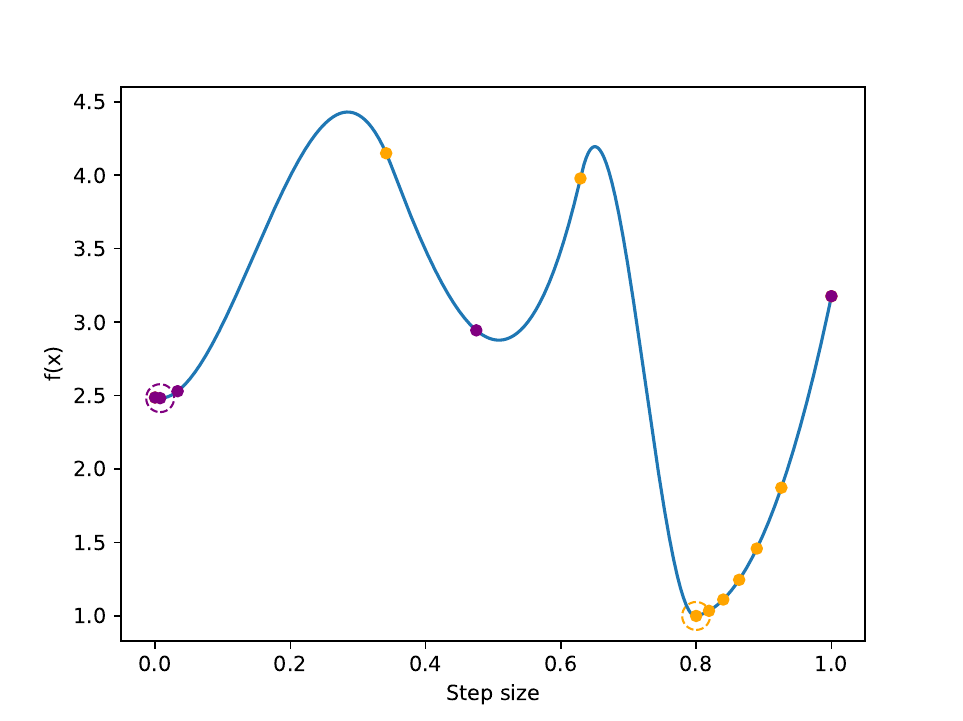}
         \caption{Example of better step length.}
         \label{fig:badplot0}
     \end{subfigure}

\hfill
     \begin{subfigure}[b]{0.45\textwidth}
         \centering
         \includegraphics[width=\linewidth]{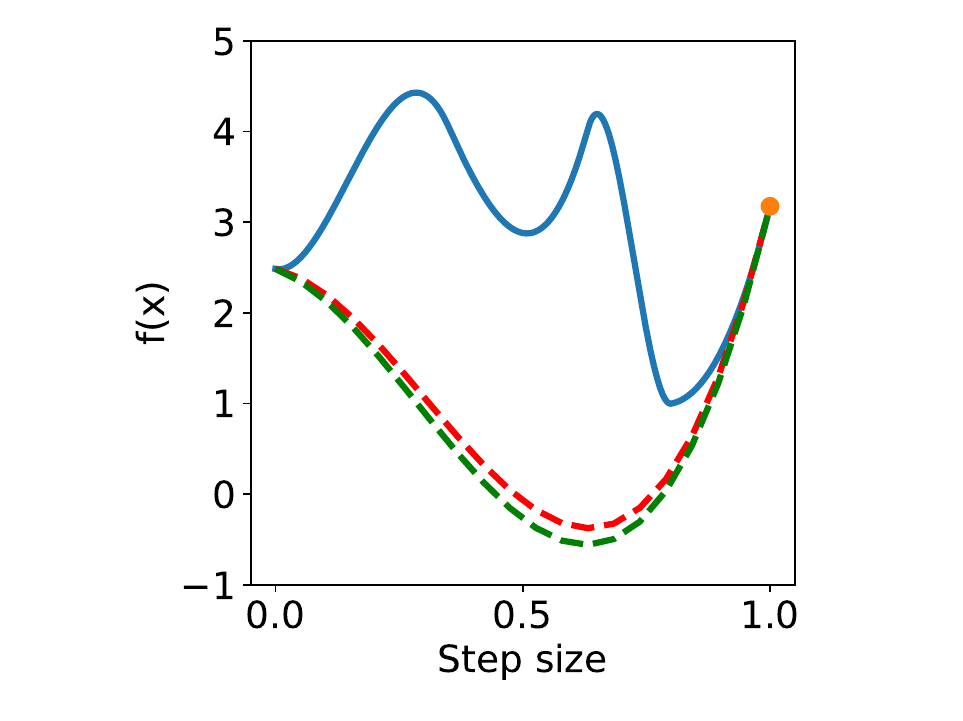}
         \caption{First iteration.}
         \label{fig:badplot1}
     \end{subfigure}
\hfill
     \begin{subfigure}[b]{0.45\textwidth}
         \centering
         \includegraphics[width=\linewidth]{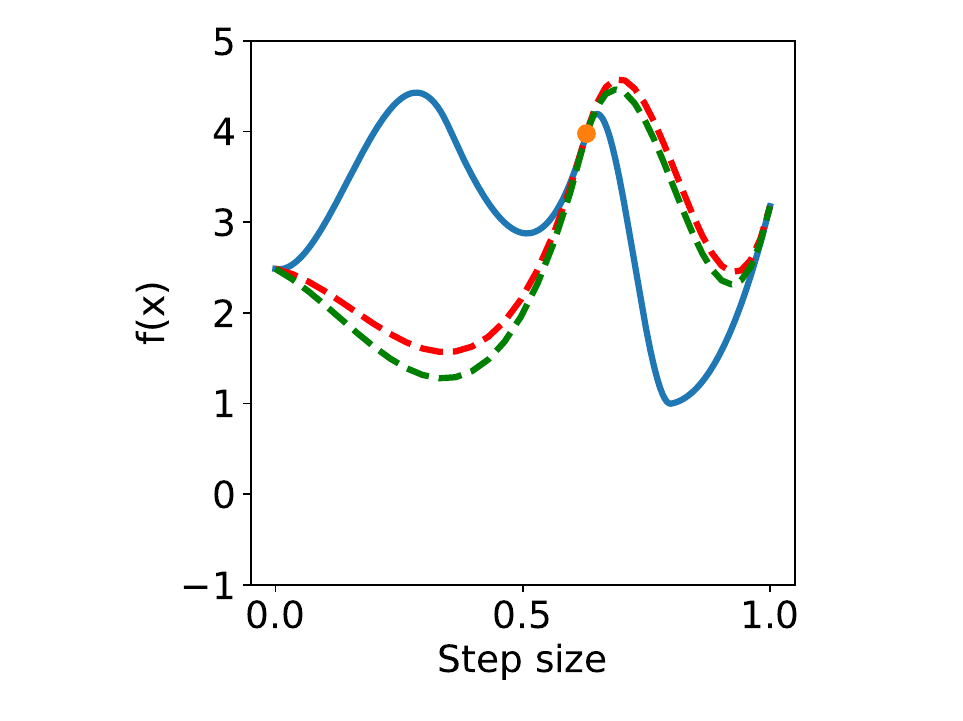}
         \caption{Second iteration.}
         \label{fig:badplot2}
     \end{subfigure}
\hfill

\hfill
     \begin{subfigure}[b]{0.45\textwidth}
         \centering
         \includegraphics[width=\linewidth]{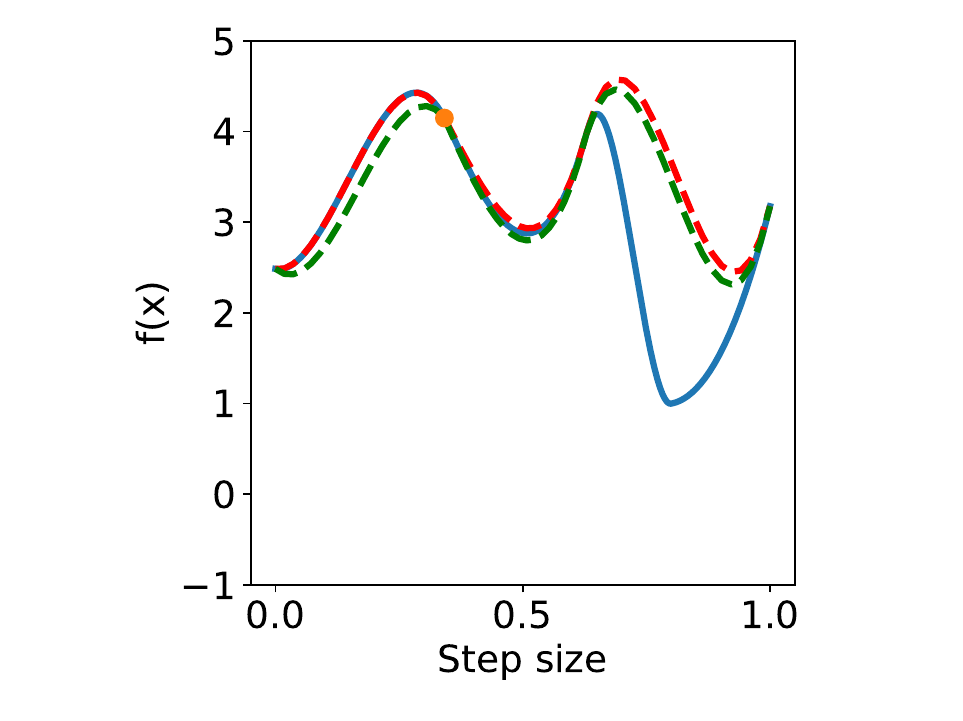}
         \caption{Third iteration.}
         \label{fig:badplot3}
     \end{subfigure}
\hfill
     \begin{subfigure}[b]{0.45\textwidth}
         \centering
         \includegraphics[width=\linewidth]{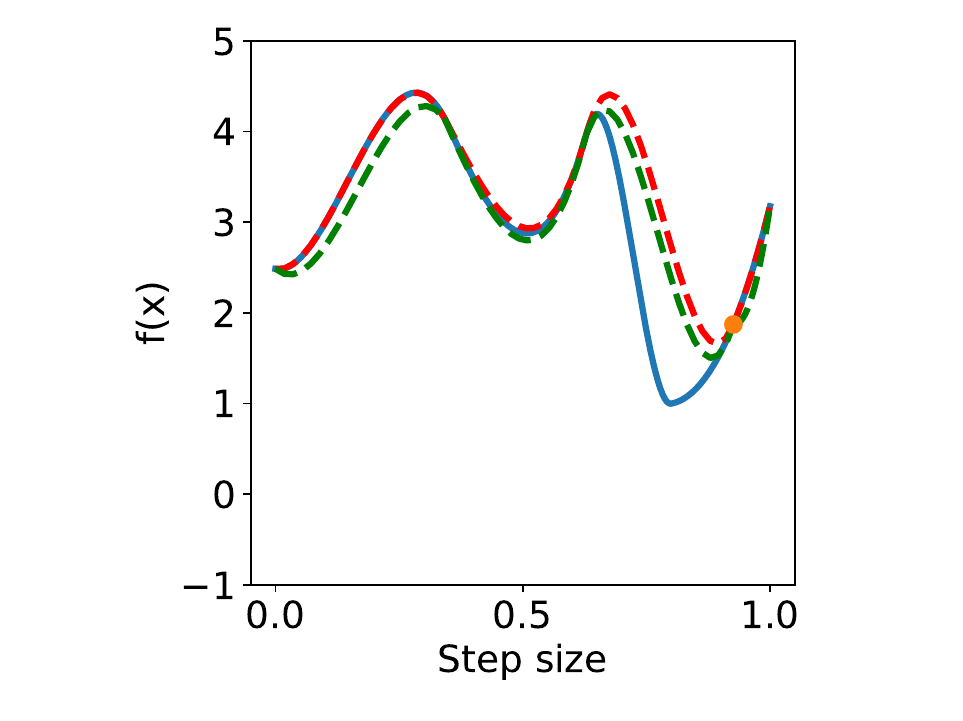}
         \caption{Fourth iteration.}
         \label{fig:badplot4}
     \end{subfigure}
\hfill
        \caption{The figure compares our line search with the Mor\'e-Thuente line search.
In Subfigure~\ref{fig:badplot0}, orange points indicate our method's query steps, while purple points indicate points tried by the Mor\'e-Thuente line search. Dashed circles indicate the final step size of the respective method.It can be seen that our method obtains a step length with a considerably better function value.
Subfigures~\ref{fig:badplot1}-~\ref{fig:badplot4} depict four iterations of our method. The blue line is the original function, the dashed red line is the cubic interpolation, and the dashed green line is the spline plus uncertainty.}
        \label{fig:2}
\end{figure}

Let us compare standard line searches with our new line search method. Standard line searches typically concentrate on a specific interval containing a point that meets the Strong Wolfe conditions, effectively narrowing their focus. This approach, however, only uses the information, i.e., function value and gradient, available from the two query points that constitute the upper and lower bounds of the interval. All the information from the other query points is discarded. Consequently, there is a risk of missing a more optimal step length. This limitation is evident when examining Fig.~\ref{fig:2}. Here, our method is shown to identify a step length that yields a significantly better function value compared to the Moré-Thuente line search. 

\newpage

\section{Experiments}
To assess the effectiveness of our proposed line search, we integrated it into the GENO solver~\cite{Geno,LaueBG22}, which originally employs a quasi-Newton method equipped with the Mor\'e-Thuente line search~\cite{MoreT94}.
We also compare against other state-of-the-art methods, including {L-BFGS-B}, which also utilizes the Mor\'e-Thuente line search, and Ipopt, which employs a filter line search~\cite{Waechter2005,Waechter2006}.
We ran all solvers on the unconstrained problems from the CUTEst test set~\cite{Gould2015}, a widely recognized benchmark featuring large-scale and numerically challenging problems.
We focus primarily on large-scale problems common in machine learning applications, where computing the Hessian is impractical or infeasible.
Consequently, all solvers could access function values and gradients but not Hessians.

We consider two criteria to evaluate convergence: the relative error in function value and the relative gradient norm.
Convergence in function value is defined as the relative error to the optimum:
\begin{equation}\label{eq:fconv}
\frac{f_{\text{solver}} - f_{\text{opt}}}{1 + |f_{\text{opt}}|} < 10^{-4}.
\end{equation}
Convergence in gradient is defined as the infinity norm of the gradient relative to the function value:
\begin{equation}\label{eq:gconv}
\frac{\|g\|_\infty}{1 + |f_{\text{solver}}|} < 10^{-6}.
\end{equation}
We declare convergence if either of these criteria is satisfied.

The experiment was conducted on 285 out of the 289 unconstrained CUTEst problems. Four problems were discarded because their optimum was unbounded. Each method was executed on each problem until convergence, with a cutoff time of 10 minutes. Convergence was determined based on the relative error in function value or the infinity norm of the gradient relative to the function value.

The results are summarized in Table~\ref{tab:solves}, which indicates the number of problems each method solved.
As can be seen in Table~\ref{tab:solves}, our proposed line search outperforms other state-of-the-art line search methods. Specifically, GENO combined with our new proposed method can solve 12 problems more than L-BFGS-B with the Mor\'e-Thuente line search, the best competing approach.

\begin{table}\caption{The table shows the number of problems each method solved.
Column {\it function conv.} shows how many problems each of the methods solved by the function value criterion (see Equation~\ref{eq:fconv}), while column {\it gradient conv.} shows how many problems each solver solved by the gradient criterion (Equation~\ref{eq:gconv}).
Column {\it convergence} shows the number of problems where either of the two equations are satisfied. }\label{tab:solves}
\centering
\begin{tabular}{lrrr}
\toprule
Solver & function conv. & gradient conv. & convergence \\
\midrule
GENO this paper & 247 & 208 & 257 \\
GENO Mor\'e-Thuente & 227 & 175 & 242 \\
Ipopt & 191 & 185 & 226 \\
{L-BFGS-B} & 231 & 211 & 245 \\
\bottomrule
\end{tabular}
\end{table}

\begin{figure}
\centering
\includegraphics[width=.6\textwidth]{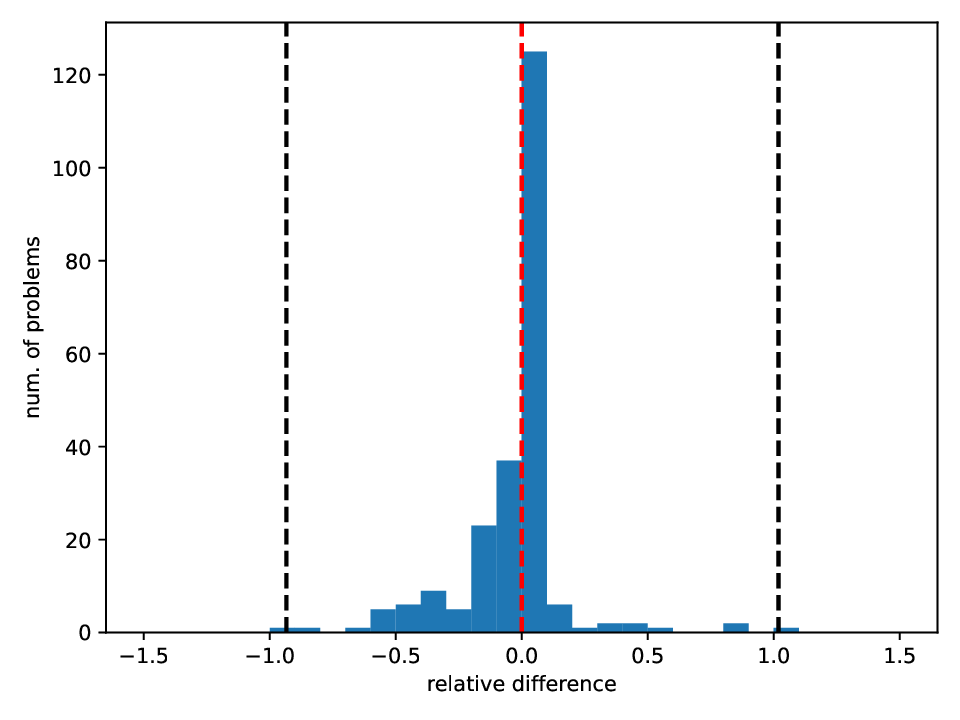}
\caption{The histogram of the relative running times of {L-BFGS-B} and GENO.
Black lines represent minimum and maximum, while the red line represents the median.
The x-axis represents the relative difference of the two methods computed as in Eqaution~\ref{eq:diff}.} \label{fig:histogram}
\end{figure}

Since {L-BFGS-B} outperformed Ipopt and GENO Mor\'e-Thuente, we conducted a detailed comparison primarily with {L-BFGS-B}. Out of the 230 problems on which both GENO and {L-BFGS-B} converged, GENO needed more function evaluations on 84 and fewer on 146 of these problems. In order to compare the number of function evaluations for one problem needed by the different solvers, we use the following formula:

\begin{equation}\label{eq:diff}
\frac{n_{\text{{L-BFGS-B}}} - n_{\text{GENO}}}{n_\text{GENO}}
\end{equation}

where \(n\) denotes the number of function evaluations performed by each solver until one of the convergence criteria was met. The distribution of relative iteration counts between GENO and {L-BFGS-B} is depicted in Figure~\ref{fig:histogram}. This figure shows that while our proposed line search, on average, uses the same number of function evaluations, it solves more problems to optimality than competing state-of-the-art approaches.

\section{Conclusion}
We presented a new line search approach. By integrating Bayesian optimization into the line search process, we have developed a method that efficiently utilizes all available data - including function values and gradients - and effectively explores the search space to identify promising regions. This approach departs from traditional line search methods, which often overlook valuable information by focusing solely on refining step length intervals. Our method's efficacy is demonstrated through rigorous testing on the CUTEst test set, outperforming existing state-of-the-art methods regarding solution quality. Our approach's simplicity and ease of implementation allow for a seamless integration into existing optimization frameworks.

\begin{credits}

\subsubsection{\discintname}
The authors have no competing interests to declare that are
relevant to the content of this article. 
\end{credits}
%
%
%

\bibliographystyle{splncs04}
\bibliography{refs}
\end{document}